\documentclass[11pt]{amsart}
\usepackage{amsmath,amssymb,amsbsy,amsfonts,amsthm,latexsym,amsopn,amstext,amsxtra,epic, euscript,amscd,indentfirst, tikz}
\usepackage{verbatim} 
\usepackage{hyperref} 
\usepackage{varioref} 
\usepackage[margin=3.5cm]{geometry}
\usepackage{pgfplots}
\usepgfplotslibrary{patchplots}
\usetikzlibrary{patterns, positioning, arrows}
\usetikzlibrary{arrows,positioning} 
\pgfdeclarepatternformonly{my crosshatch dots}{\pgfqpoint{-1pt}{-1pt}}{\pgfqpoint{5pt}{5pt}}{\pgfqpoint{6pt}{6pt}}%
{
    \pgfpathcircle{\pgfqpoint{0pt}{0pt}}{.5pt}
    \pgfpathcircle{\pgfqpoint{3pt}{3pt}}{.5pt}
    \pgfusepath{fill}
}

\DeclareFontFamily{U}{mathc}{}
\DeclareFontShape{U}{mathc}{m}{it}%
{<->s*[1.03] mathc10}{}

\DeclareMathAlphabet{\mathscr}{U}{mathc}{m}{it}

\begin{document}

\newtheorem{theorem}{Theorem}
\newtheorem*{theorem*}{Theorem}
\newtheorem{conjecture}[theorem]{Conjecture}
\newtheorem{proposition}[theorem]{Proposition}
\newtheorem*{proposition*}{Proposition}
\newtheorem{question}[theorem]{Question}
\newtheorem{lemma}[theorem]{Lemma}
\newtheorem{remark}[theorem]{Remark}
\newtheorem{definition}[theorem]{Definition}
\newtheorem{example}[theorem]{Example}
\newtheorem*{definition*}{Definition}
\newtheorem{cor}[theorem]{Corollary}
\newtheorem*{cor*}{Corollary}
\newtheorem*{result*}{Result}
\newtheorem{obs}[theorem]{Observation}
\newtheorem{proc}[theorem]{Procedure}
\newcommand{\comments}[1]{}
\newcommand{\todo}[1]{\textbf{\textcolor{red}{[#1]}}}

\def\Z{\mathbb Z}
\def\Za{\mathbb Z^\ast}
\def\Fq{{\mathbb F}_q}
\def\R{\mathbb R}
\def\N{\mathbb N}
\def\i{\sqrt{-1}}
\def\k{\kappa}
\def\M{\mathcal{M}}

\title[Optimal Transport and Complex Geometry]{The K\"ahler Geometry of Certain Optimal Transport Problems}

 \author[University of Michigan]{Gabriel Khan} 
 \author[University of Michigan]{Jun Zhang}
 
\email{gabekhan@umich.edu}
\email{junz@umich.edu}

\date{\today}

\maketitle

\begin{abstract}
Let $X$ and $Y$ be domains of $\mathbb{R}^n$ equipped with probability measures $\mu$ and $ \nu$, respectively. We consider the problem of optimal transport from $\mu$ to $\nu$ with respect to a cost function $c: X \times Y \to \mathbb{R}$. To ensure that the solution to this problem is smooth, it is necessary to make several assumptions about the structure of the domains and the cost function. In particular, Ma, Trudinger, and Wang \cite{MTW} established regularity estimates when the domains are strongly \textit{relatively c-convex} with respect to each other and the cost function has non-negative \textit{MTW tensor}. For cost functions of the form $c(x,y)= \Psi(x-y)$ for some convex function $\Psi$, we find an associated K\"ahler manifold whose orthogonal anti-bisectional curvature is proportional to the MTW tensor. We also show that relative $c$-convexity geometrically corresponds to geodesic convexity with respect to a dual affine connection. Taken together, these results provide a geometric framework for optimal transport which is complementary to the pseudo-Riemannian theory of Kim and McCann \cite{KM}.

We provide several applications of this work. In particular, we find a complete K\"ahler surface with non-negative orthogonal anti-bisectional curvature that is not a Hermitian symmetric space or biholomorphic to $\mathbb{C}^2$. We also address a question in mathematical finance raised by Pal and Wong \cite{MSP} on the regularity of \textit{pseudo-arbitrages}, or investment strategies which outperform the market.


\end{abstract}

\section{Introduction}

Optimal transport is a classic field of mathematics combining ideas from geometry, probability, and analysis. The problem was first formalized by Gaspard Monge in 1781 \cite{Monge}. In his work, he considered a worker who is tasked with moving a large pile of sand into a prescribed configuration and wants to minimize the total effort required to complete the job. Trying to determine the optimal way of transporting of the sand leads into deep and subtle mathematical phenomena and is a thriving field of research. Furthermore, optimal transport has many practical applications. Monge's work was originally inspired by a problem in engineering, but these same ideas can be applied to logistics, economics, computer imaging processing, and many other fields \cite{PC}.

The modern framework for optimal transport, due to Kantorovich \cite{LVK}, considers arbitrary couplings between two probability measures. In this formulation, we consider $X$ and $Y$ as Borel subsets of two metric spaces equipped with probability measures $\mu$ and $\nu$, respectively. Intuitively, $d \mu$ is the shape of the original sand pile and $d \nu$ is the target configuration. To transport the sand from $\mu$ to $\nu$, we consider a \textit{coupling} of $\mu$ and $\nu$, which is a non-negative measure on $X \times Y$ whose marginal distributions are $\mu$ and $\nu$, respectively.
To measure the efficiency of a plan for transport $\mu$ to $\nu$, we consider a lower-semicontinuous cost function $c: X \times Y \to \mathbb{R}$. The solution to the Kantorovich optimal transport problem is the coupling $\gamma$ which achieves the smallest total cost
\[ \min_{\gamma \in \Gamma(\mu,\nu)} \int_{X \times Y} c(x,y) d \gamma(x,y). \] 
Here $\Gamma(\mu,\nu)$ is the set of all couplings of $\mu$ and $\nu$, (i.e. joint probabilities whose marginal distributions are $\mu$ and $\nu$, respectively). In this case, a minimizing measure $\gamma$ is referred to as the \textit{optimal coupling}. An optimal coupling exists for very general measures and cost functions, so the Kantorovich approach is a flexible and powerful framework to study optimal transport

In Monge's work, it is assumed that the mass at a given point will not be subdivided and sent to multiple locations. This is known as \textit{deterministic} optimal transport, and seeks to find a measurable map $\mathbb{T}:X \to Y$ so that the optimal coupling is entirely contained within the graph of $\mathbb{T}$. When this occurs, the map $\mathbb{T}$ is known as the \textit{optimal map}.  A priori, there is no guarantee that optimal transport is deterministic, so a Monge solution may not exist for a given optimal transport problem. However, we will discuss certain sufficient conditions for the optimal transport to be deterministic in Section 2. 

For deterministic optimal transport, it is natural to ask whether the optimal map is continuous or even smooth. This is known as the \textit{regularity problem for optimal transport}. Historically, most of the work on this problem was done in Euclidean space for the cost $c(x,y) = \|x-y\|^p$  better known as the $p$-Wasserstein cost.
 
For more general cost functions (such as Wasserstein costs on Riemannian manifolds), the groundbreaking work was done by Ma, Trudinger and Wang \cite{MTW}, who proved that the transport map is smooth under the assumptions that a certain non-linear fourth order quantity, known as the MTW tensor (denote $\mathfrak{S}$), is non-negative and that $X$ and $Y$ are \textit{relatively c-convex} with respect to each other. These results were refined by Loeper \cite{GL}, who showed that the non-negativity of $\mathfrak{S}$ is necessary to even establish continuity for smooth measures $\mu$ and $\nu$. Furthermore, he gave some insight into the geometric significance of the MTW tensor. Later work of Kim and McCann \cite{KM} furthered this understanding by presenting a pseudo-Riemannian framework for optimal transport in which the MTW tensor is the curvature of certain light-like planes. 

\subsection{Our Results}

 In this paper, we mainly consider $\Psi$-costs, which we define as follows.
\begin{definition*}[$\Psi$-cost]\label{Psi costs}
Let $\Psi: \M \to \mathbb{R}$ be a convex function on an open domain $\M$ in Euclidean space. For open domains $X$ and $Y$ in $\mathbb{R}^n$, a $\Psi$-cost is a cost function of the form 
\begin{eqnarray*}
c:  X \times Y & \to & \mathbb{R} \\
  c(x,y) & = & \Psi(x-y) 
\end{eqnarray*}
\end{definition*}

 These costs were previously studied by Gangbo and McCann \cite{GM} and by Ma, Trudinger and Wang \cite{MTW}. For such a cost to be well defined, $\M$ must contain the difference set $X-Y$, defined as
\[X-Y := \{ z \in \mathbb{R}^n ~|~ \exists \, x \in X,~y \in Y \textrm{ such that } z=x-y \}. \]

We can now summarize the main results of our work, which associate a complex manifold to a given $\Psi$-cost. To do so, we consider $\M$ as a Hessian manifold, using $\Psi$ as its potential function. Such manifolds naturally admit a dual pair of flat connections, which we denote $D$ and $D^*$. Using the primal flat connection $D$, the tangent bundle
$T \M$ can be equipped with a K\"ahler metric, known as the Sasaki metric and denoted $(T \M, g^D, J^D)$. Our main result shows the following correspondence between the curvature of this metric and the MTW tensor.

\begin{theorem*}
Let $X$ and $Y$ be open sets in $\mathbb{R}^n$ and $c$ be a $\Psi$-cost. Then the MTW tensor $\mathfrak{S}$ of $c$ satisfies the following identity:
\[ \frac{1}{2}\mathfrak{S}(\eta,\xi)=  \mathfrak{R}_{g^D}(\xi, J \eta^\sharp , \xi, J \eta^\sharp)-  \mathfrak{R}_{g^D}(\eta^\sharp, \xi, \eta^\sharp, \xi)  \]
where $\xi$ and $\eta$ are an orthogonal real vector-covector pair (which we extend to $T\M$) and $\mathfrak{R}_{g^D}$ is the curvature of $(T \M, g^D, J^D)$.
\end{theorem*}

 For reasons that we will explain later, we call the right hand expression the \textit{orthogonal anti-bisectional curvature}. We furthermore show that relative $c$-convexity of sets is geodesic convexity with respect to the dual affine connection on $\M$.
 
 \begin{proposition*}
 For a $\Psi$-cost, a set $Y$ is $c$-convex relative to $X$ if and only if, for all $x \in X$, the set $x-Y$ is geodesically convex with respect to the dual connection $D^*$. 
 \end{proposition*}
 
 Apart from providing a new geometric framework for the regularity problem, we can use these results to address several questions of independent interest.

\subsubsection{Applications to Complex Geometry}
This approach can be used to construct several examples of interesting metrics with subtle non-negativity properties. In particular, we find a complete complex surface which is neither biholomorphic to $\mathbb{C}^2$ nor Hermitian symmetric but whose orthogonal anti-bisectional curvature is non-negative. Many of the complex manifolds constructed using this approach are of independent interest, and we will provide a few examples which we will study in depth in future work.

\subsubsection{Applications to Mathematical Finance}
Our second main application is to establish regularity for a certain problem in portfolio design theory. Recent work of Pal and Wong \cite{GRA} studies the problem of finding \textit{pseudo-arbitrages}, which are investment strategies that outperform the market almost surely in the long run. Their work shows that this is equivalent to solving an optimal transport problem with a statistical divergence that is closely related to the free energy in statistical physics.

For this problem, our approach relates the MTW tensor of this cost to a K\"ahler manifold with constant positive holomorphic sectional curvature. As such, this cost function satisfies the $MTW(0)$ condition (and also satisfies a stronger condition known as non-negative cost-curvature). We further show that relative $c$-convexity corresponds precisely to the standard notion of convexity on the probability simplex. Combining these calculations, we can apply the results of \cite{TW} to obtain a regularity theory of portfolio maps and their associated displacement interpolations. This addresses a question asked in \cite{MSP}, and intuitively shows that when the market conditions change slightly, the investment strategy similarly does not change by much.

\subsubsection{$\mathcal{D}_\Psi^{(\alpha)}$-Divergences and Information Geometry}

Although our main results are stated in terms of $\Psi$-costs, they also hold (with minor modifications) for cost functions that are $\mathcal{D}_\Psi^{(\alpha)}$-divergences, which were previously studied by the second author \cite{JZ}. 

\begin{definition*}[$\mathcal{D}_\Psi^{(\alpha)}$-divergence]\label{Dalphadiv}
Let $\Psi: \M \to \mathbb{R}$ be a convex function on a convex domain $\M$ in Euclidean space. For two points $x, y \in \M$ and  $\alpha \in (-1,1)$,  a $\mathcal{D}_\Psi^{(\alpha)}$-divergence is a function of the form
\begin{eqnarray*}
\mathcal{D}_{\Psi}^{(\alpha)}(x , y)  &= & \frac{4}{1-\alpha^2} \left[ \frac{1-\alpha}{2} \Psi(x) +  \frac{1+\alpha}{2} \Psi( y) - \Psi \left( \frac{1-\alpha}{2} x+ \frac{1+\alpha}{2}  y \right) \right].  
\end{eqnarray*}
\end{definition*}

As $\alpha$ approaches $\pm 1$, the $\mathcal{D}_\Psi^{(\alpha)}$-divergence converges to the Bregman divergence \cite{LB}. When $\alpha =0$ and $\Psi$ is quadratic, this divergence simply becomes the 2-Wasserstein distance. More broadly, divergences are a generalization of distance functions, where the assumptions of symmetry and the triangle inequality are dropped. They are widely used in statistics and information geometry, which studies the geometry of parametrized families of probability distributions. Although a background on the subject is not strictly necessary for this paper, much of the motivation for this work is information geometric, and we will freely use basic results about exponential families to provide intuition and context for our applications. For a more complete background on information geometry, we refer readers to the book by Amari \cite{SIA}.


Both $\Psi$-costs and $\mathcal{D}_\Psi^{(\alpha)}$-divergences involve a convex function defined on an open domain of $\mathbb{R}^n$. The primary difference is whether it is necessary to assume that $X-Y \subset \M$ (as for a $\Psi$-cost), or that $X,\,Y,$  and $\frac{1-\alpha}{2} X+ \frac{1+\alpha}{2}  Y \subset \M$ (as for a $\mathcal{D}_\Psi^{(\alpha)}$-divergence). In order to ensure that the $\mathcal{D}_\Psi^{(\alpha)}$-divergence is well defined, we will assume that the domain of $\Psi$ is convex. For many of the relevant examples, $\Psi$ will be the so-called \textit{log-partition function} of an exponential family in its natural parameters. It is a general property of exponential families that the domains of such functions are convex, so this convexity assumption will be automatically satisfied. As such, the geometry of the $\mathcal{D}_\Psi^{(\alpha)}$-divergences is more natural, as there is no need to consider the difference set $X-Y$. A preliminary announcement of our results (stated in terms of $\mathcal{D}_\Psi^{(\alpha)}$-divergences) has been
accepted to the 2019 proceedings of the Geometric Science of Information \cite{KZ}.





\subsection{Layout of the Paper}

In Section 2 we discuss some background information on optimal transport. Section 3 discusses some background information on Hessian manifolds and the curvature of the Sasaki metric. Both of these sections are largely review and can be skipped by readers familiar with the theory. In Section 4, we state our main results, which show the precise interaction between complex/information geometry and the regularity theory of optimal transport. In Section 5, we explore various applications of this result. In Section 6, we conclude with a section of open questions, which we hope to explore in future work.

\subsection{Notation}

We have attempted to preserve the notation from \cite{DPF} and \cite{HS} as much as possible, while minimizing abuse of notation or overlap. For clarity, we introduce some notational conventions now.

Throughout the paper, $X$ and $Y$ will denote open domains in $\mathbb{R}^n$. Invariably, these will be smooth and bounded. We will use $ \{ x^i \}_{i=1}^n$ as coordinates on $X$ and $\{ y^i \}_{i=1}^n$ as coordinates on $Y$. To study optimal transport, we will use $c(x,y)$ to denote a lower-continuous cost function $c:X \times Y \to \R$. Oftentimes, the domain of $c$ will be larger than $X \times Y$, but we will often ignore this. To avoid confusion with coordinate functions and the notation for tangent spaces, we denote the solutions to Monge-Ampere type equations as $\textsc{u}$, and the associated optimal map $\mathbb{T}_\textsc{u}$.

For the most part, $\M$ will be an open domain in Euclidean space, and $\Psi$ will denote a convex function $\Psi: \M \to \R$.  
It is instructive to also consider $\M$ as a manifold, and we will use $\{ u^i \}_{i=1}^n$ as its coordinates. In Section 5, we will occasionally need to square the coordinate functions. When doing so, we denote coordinate functions with subscripts (i.e. $\{ u_i \}_{i=1}^n$). When considering the tangent bundle of $\M$ (denoted $T \M$), we will use bundle coordinates $\{ (u^i,v^i) \}_{i=1}^n$. This notation is a change from Satoh \cite{HS}, and is done to avoid overusing `$x$' and `$y$'.

In order to prescribe $T \M$ with an almost Hermitian structure, it is necessary to consider an affine connection on $\M$, which we denote by $D$. Furthermore, we use $\mathcal{W}, \mathcal{X}, \mathcal{Y}, \mathcal{Z} $ and $\xi$ to denote tangent vectors on $\M$ (i.e. elements of $T \M$). This is the convention of \cite{HS}, except with calligraphic font to avoid confusion with our notation for domains. When computing the MTW tensor, we will denote the vectors in the MTW tensor by $\xi$ and the covectors by $\eta$.


To simplify the derivative notation, for a two variable function $c(x,y)$, we use $c_{I,J}$ to denote 
$\partial_{x^I} \partial_{y^J} c$ for multi-indices $I$ and $J$. Furthermore, $c^{i,j}$ denotes the matrix inverse of the mixed derivative $c_{i,j}$. For a convex function $\Psi$, we use the notation $\Psi_J$ to denote 
$\partial_{u^J} \Psi$ for a multi-index $J$ and the notation $\Psi^{ij}$ to denote the matrix inverse of $\Psi_{ij}$. Finally, we will use Einstein summation notation throughout the paper.

\subsection{Acknowledgments} 

The first author would like to thank Fangyang Zheng and Bo Guan for their helpful discussions with this project. The project is supported by DARPA/ARO Grant W911NF-16-1-0383 (“Information Geometry: Geometrization of Science of Information”, PI: Zhang). A preliminary announcement of these results has been accepted to the proceedings of the Geometric Science of Information 2019 \cite{KZ}.

\section{Background on the regularity theory of optimal transport}
The main focus of our paper is to study the assumptions needed to ensure optimal transport is regular. In order to understand these, we first review several preliminary results on the regularity theory of optimal transport.

As our primary interest is the geometric structure of the regularity problem, we will not make use of the sharpest possible regularity estimates. The material in this section is based off the survey paper by De Phillipis and Figalli \cite{DPF}, which gives a more complete overview to the regularity theory. For a more complete overview on optimal transport, see the book by Villani \cite{OTON}.


Recall that the regularity problem is of interest primarily in the case of deterministic
optimal transport. As such, we first discuss some conditions that ensure the Kantorovich
optimal transport problem is deterministic. The following theorem was originally proven by Brenier for the 2-Wasserstein cost \cite{YB} and proved in more generality by Gangbo and McCann \cite{GOT}. It gives sufficient conditions for deterministic transport and shows that the optimal maps can be found by solving a Monge-Ampere type equation.

\begin{theorem}
Let $X$ and $Y$ be two open subsets of $\R^n$ and consider a cost function $c:X \times Y \to \R$. Suppose that $ d \mu$ is a smooth probability density supported on $X$ and that $d \nu $ is a smooth probability density supported on $Y$. Suppose that the following conditions hold:

\begin{enumerate}
\item The cost function $c$ is of class $C^4$ with $\| c \|_{C^4(X \times Y)} < \infty$
\item For any $x \in X$, the map $ Y \ni y \to  c_x(x,y) \in \R^n$ is injective.
\item For any $y \in Y$, the map $ X \ni x \to  c_y(x,y) \in \R^n$ is injective.
\item $\det(c_{x,y})(x,y) \neq 0$ for all $(x,y) \in X \times Y$.
\end{enumerate}

Then there exists a $c$-convex function $\textsc{u}: X \to \R$ such that the map $\mathbb{T}_\textsc{u} : X \to Y$ defined by $\mathbb{T}_\textsc{u}(x) := c\textrm{-}\exp_x(\nabla \textsc{u}(x))$ is the unique optimal transport map sending $\mu$ onto $\nu$. Furthermore, $\mathbb{T}_\textsc{u}$ is injective $d \mu$-a.e.,
\begin{equation} \label{Monge Ampere}
 | \det(\nabla \mathbb{T}_\textsc{u}(x))| = \frac{ d \mu (x)}{ d \nu(\mathbb{T}_\textsc{u}(x))} \hspace{.2in} d \mu-a.e., 
\end{equation}
and its inverse is given by the optimal transport map sending $\nu$ onto $\mu$.

\end{theorem}

In order to express Equation \ref{Monge Ampere} more concretely, we define the $c$-exponential map (denoted $c\textrm{-}\exp_x$).

\begin{definition*} [$c$-exponential map]
For any $x \in X, y \in Y, p \in \mathbb{R}^n$, the $c$-exponential map satisfies the following identity.
\[ c\textrm{-}\exp_x(p) = y \iff p = - c_x(x, y). \]
\end{definition*}

 For the 2-Wasserstein cost on a Riemannian manifold, the $c$-exponential is exactly the standard exponential map, which motivates its name. For this cost in Euclidean space, Equation \ref{Monge Ampere} becomes the standard Monge-Ampere equation
\[ \det \left( \nabla^2 \textsc{u}(x) \right) = \frac{f(x)}{g(\nabla \textsc{u} (x))}. \] 
Due to this simple form for Equation \ref{Monge Ampere}, much of the initial work on the regularity problem was done for the squared-distance cost in Euclidean space. In this setting, Caffarelli \cite{LC} and others proved a priori estimates under certain convexity and smoothness assumptions on the measures (for a more complete history, see \cite{DPF}). Caffarelli also observed there is no hope of proving interior regularity without assuming that the support of the target measure is convex.

 For more general cost functions, Ma, Trudinger and Wang's breakthrough work in 2005 \cite{MTW} gave three conditions that ensure $C^2$ regularity for the solutions of Monge-Ampere equations. In this paper, we will use a stronger version of this result, originally proven by Trudinger and Wang \cite{TW}.

\begin{theorem} \label{MTW}
Suppose that $c:X \times Y \to \R$, $\mu$, and $\nu$ satisfy the hypothesis of Theorem 1, and that the densities $d \mu$ and $ d \nu$ are bounded away from zero and infinity on their respective supports $X$ and $Y$. Suppose further that the following holds.

\begin{enumerate}
\item $X$ and $Y$ are smooth.
\item The domain $X$ is strictly $c$-convex relative to the domain $Y$. 
\item The domain $Y$ is strictly $c^*$-convex relative to the domain $X$. 
\item 
For all vectors $\xi, \eta \in \R^n$ with $\xi \perp \eta$, the following inequality holds.
\begin{equation}\label{MTW0}
\mathfrak{S}(\xi,\eta) := \sum_{i,j,k,l,p,q,r,s} (c_{ij,p}c^{p,q}c_{q,rs}-c_{ij,rs})c^{r,k}c^{s,l} \xi^i \xi^j \eta^k \eta^l \geq 0 
\end{equation}
\end{enumerate}

Then $\textsc{u} \in C^\infty(\overline X)$ and $\mathbb{T}_\textsc{u}: \overline{X} \to \overline{Y}$ is a smooth diffeomorphism, where $\mathbb{T}_\textsc{u}(x) = c\textrm{-}\exp_x(\nabla \textsc{u}(x))$.
\end{theorem}

While we will not discuss the proof in detail, we note that the main challenge is obtaining an a priori $C^2$ estimate on \textsc{u}. Once such an estimate is established, the Monge-Ampere equation can be linearized at \textsc{u}, at which point standard elliptic bootstrapping yields estimates of all orders and implies that $\mathbb{T}_\textsc{u}$ is smooth. 

The main results of this paper study the assumptions of Theorem \ref{MTW}, so we discuss these in more detail. The first condition is self-explanatory, while the second and third provide the proper notions of convexity for the supports. To explain this in detail, we define the notion of $c$-convexity for sets.

\begin{definition*} [$c$-segment] A $c$-segment in $X$ with respect to a point $y$ is a solution set $\{x\}$ to $c_y(x, y) \in \ell$ for $\ell$ a line segment in $\R^n$. A $c^*$-segment in $Y$ with respect to a point $x$ is a solution set $\{ y \}$ to $c_x(x, y) \in \ell$ where $\ell$ is a line segment in $\mathbb{R}^n$.
\end{definition*}

\begin{definition*} [$c$-convexity]
A set $E$ is $c$-convex relative to a set $E^*$ if for any two points $x_0, x_1 \in E$ and
any $y \in E^*$, the $c$-segment relative to $y$ connecting $x_0$ and $x_1$ lies in E. Similarly we say $E^*$ is
$c^*$-convex relative to $E$ if for any two points $y_0, y_1 \in E^*$ and any $x \in E$, the $c^*$-segment relative to
$x$ connecting $y_0$ and $y_1$ lies in $E^*$.
\end{definition*}


Finally, we discuss Inequality \ref{MTW0}, which is known as the $MTW(0)$ condition and is a weakened version of the $MTW(\kappa)$ condition.

\begin{definition*}[$MTW(\kappa)$]
A cost function $c$ satisfies the $MTW(\kappa)$ condition if for  any orthogonal vector-covector pair $\eta$ and $\xi$, $\mathfrak{S}(\eta,\xi) \geq  \kappa|\eta|^2 |\xi|^2 $ for $\kappa>0$.
\end{definition*}


 Ma, Trudinger and Wang's original work relied on the $MTW(\kappa)$ assumption, and this stronger condition is used in many applications. Although it is not immediately apparent, $\mathfrak{S}(\xi,\eta)$ is tensorial (coordinate-invariant) so long as one considers $\eta$ as a covector \cite{KM}, which we will do throughout the rest of the paper. Furthermore, it transforms quadratically in $\eta$ and $\xi$, but is highly non-linear and non-local in the cost function. 

The geometric significance of the MTW tensor is an active topic of research. On a Riemannian manifold, Loeper \cite{GL} gave some insight into its behavior. His work showed that for the 2-Wasserstein cost, the MTW tensor is proportional to the sectional curvature on the diagonal $x=y$. In this paper, he also showed that $c$-convexity and non-negativity of the MTW tensor are essentially necessary conditions to prove regularity of optimal transport.

Building off of Loeper's results, Kim and McCann gave a geometric framework for optimal transport \cite{KM}. In their formulation, optimal transport is expressed in terms of a pseudo-Riemannian metric on the manifold $X \times Y$ and the MTW tensor becomes the curvature of light-like planes. This interpretation holds for arbitrary cost functions, which gives intrinsic geometric structure to the regularity problem. Our geometric interpretation is different, but many of the formulas appear similar, in part due to the fact that Kim and McCann chose notation reminiscent of complex geometry. 

Before concluding our background discussion on optimal transport, we will introduce one more strengthening of the $MTW(0)$ condition, known as non-negative ``cross-curvature" \cite{FKM}.

\begin{definition*}[Non-negative cross-curvature]
A cost function $c$ has non-negative (resp. strictly positive) cross-curvature if, for any vector-covector pair $\eta$ and $\xi$, 
\[ \mathfrak{S}(\eta,\xi) \geq  0 \textrm{ (resp. } \kappa |\eta|^2 |\xi|^2). \]
\end{definition*}

Note that non-negative cross-curvature is stronger than $MTW(0)$, as the non-negativity must hold for all pairs $\eta$ and $\xi$, not simply orthogonal ones. Cross-curvature was introduced by Figalli, Kim, and McCann \cite{FKM} to study a problem in microeconomics. In later work, they also showed that stronger regularity for optimal maps can be proven with this assumption \cite{FKM2}. Cross-curvature was also studied by Sei \cite{TS} for an application in statistics.



\section{Hessian manifolds and the Sasaki metric}

In order to interpret the MTW tensor as a complex-geometric curvature, we must
study the Sasaki metric, which is an almost-Hermitian structure on the tangent bundle
of a Riemannian manifold. We now discuss some background on this metric, focusing
on the case of Hessian manifolds, for which the Sasaki metric is K\"ahler.

\subsection{The Sasaki Metric on the Tangent Bundle}

On a general Riemannian manifold $(\M,g)$ with an affine connection $ D$, the tangent bundle naturally inherits an almost-Hermitian structure $(T \M, g^{D}, J^{D})$. This is known as the \textit{Sasaki metric} and was introduced by Dombrowski \cite{PD}. For completeness, we will present a brief overview of this construction, derived from the work of Satoh \cite{HS}.

 Given local coordinates $\{u^i\}_{i=1}^n$ on $\M$, we denote the Christoffel symbols of the connection by $\Gamma_{ji}^k$ where
\[ D_{\partial_{u^i}} \partial_{u^j} := \Gamma_{ji}^k \partial_{u^k}. \]

On the tangent bundle $T\M$, we can define smooth functions $v^1, \ldots, v^n$ by $v^j(\mathcal{X}) = \mathcal{X}^j$ for a vector $\mathcal{X}=\mathcal{X}^i \partial_{u^i}$. The collection of functions $\{(u^i,v^i) \}$ form local coordinates for $T\M$. 

For a tangent vector $\xi \in T_u \M$ (which we consider as a point in the tangent bundle $T \M$) and a tangent vector $\mathcal{X}= \mathcal{X}^i \partial_{u^i} \in T_u \M$, we can define vertical and horizontal lifts of $\mathcal{X}$ at $\xi$, denoted $\mathcal{X}_\xi^V$ and $\mathcal{X}_\xi^H$, respectively. These are elements of $T_\xi(T\M)$, which are defined as follows:
\[\mathcal{X}_\xi^H = \mathcal{X}^i \partial_{u^i} - \Gamma_{ij}^k \mathcal{X}^i v^j(\xi) \partial_{v^k}, \]
\[ \mathcal{X} _ { \xi } ^ { V } =\mathcal{X} ^ { i }  \partial_ {v ^ { i } }. \]

This yields a decomposition of $T_\xi(T\M)$ into horizontal and vertical subspaces, which depends on the choice of connection $D$:
 \[ T_\xi(T\M)= H_\xi (T\M) \oplus V_\xi(T\M). \]
As such, there is a natural identification  $H_\xi (T\M) \cong V_\xi(T\M) \cong T_u \M$, which we use to construct the Sasaki metric  (Definition 2.1 of \cite{HS}).

\begin{definition*}[Sasaki metric] \label{Sasaki metric}
Let $(\M^n,g)$ be a Riemannian manifold with affine connection $D$. The Sasaki metric is the following almost-Hermitian structure on $T \M$.

For $\mathcal{X},\mathcal{Y} \in T_u \M$ and $\xi \in T \M$ with $\xi = (u,v)$ in bundle coordinates, the almost complex structure $J^D$ is defined as 
\[ J ^ { D } \mathcal{X} _ { \xi } ^ { H } = \mathcal{X} _ { \xi } ^ { V } , \quad J ^ { D } \mathcal{X} _ { \xi } ^ { V } = - \mathcal{X} _ { \xi } ^ { H } \]
and the Riemannian metric $\tilde g^D$ is defined as
\[\widetilde { g } ^ { D } \left( \mathcal{X} _ { \xi } ^ { H } , \mathcal{Y} _ { \xi } ^ { H } \right) = \widetilde { g } ^ { D } \left( \mathcal{X} _ { \xi } ^ { V } , \mathcal{Y} _ { \xi } ^ { V } \right) = g ( \mathcal{X} , \mathcal{Y} ) , \quad \widetilde { g } ^ { D } \left( \mathcal{X} _ { \xi } ^ { H } , \mathcal{Y} _ { \xi } ^ { V } \right) = 0.\]
\end{definition*}

This induces an almost-Hermitian structure on $T\M$, which depends on both the choice of metric and connection on $\M$. A priori, this structure is neither integrable nor almost K\"ahler. Work by Dombrowski and by Satoh give sufficient and necessary conditions for these properties to hold.

\begin{theorem}[\cite{PD} \cite{HS}]

Let $(\M, g)$ be a Riemannian manifold with an affine connection $D$. The almost-Hermitian manifold $(T \M, g^D, J^D)$ satisfies the following.
\begin{enumerate}
\item The almost-Hermitian structure is integrable if and only if the connection $D$ is flat \cite{PD}.
\item The almost-Hermitian structure is almost-K\"ahler if and only if the dual connection $D^*$ is torsion free \cite{HS}.
\item The almost-Hermitian structure is K\"ahler if and only if  $(D, g)$ are dually flat, which further implies that $g$ is Hessian \cite{PD}.
\end{enumerate}

\end{theorem}

\subsection{Hessian Manifolds}

We are primarily interested in the case where $T\M$ is K\"ahler, for which we must study Hessian manifolds (also known as \textit{affine-K\"ahler} manifolds). There are two equivalent definitions for such manifolds; with the former definition primarily used in differential geometry and the latter primarily used in information geometry. 

\begin{definition*}[Differential geometric]
A Riemannian manifold $(\M,g)$ is Hessian if there is an atlas of local coordinates $\{u^i \}_{i=1}^n$ so that for each coordinate chart, there is a convex potential $\Psi$ such that
$$g_{ij}= \frac{ \partial^2 \Psi}{\partial u^i u^j}.  $$
Furthermore, the transition maps between these coordinate charts are affine (i.e. $\M$ is an affine manifold).
\end{definition*}

\begin{definition*}[Information geometric]
 A Riemannian manifold $(\M,g)$ is said to be \textit{Hessian} if it admits dually flat connections. That is to say, it admits two flat (torsion- and curvature-free) connections $D$ and $D^*$ satisfying
\begin{equation} \label{conjugateconnection}
\mathcal{X}(g(\mathcal{Y},\mathcal{Z})) = g(D_\mathcal{X} \mathcal{Y}, \mathcal{Z}) + g(\mathcal{Y}, D^*_\mathcal{X}  \mathcal{Z}) 
\end{equation}
 for all vector fields $\mathcal{X},~ \mathcal{Y}$, and  $\mathcal{Z}$.
\end{definition*}

 Although these definitions initially appear different, they are actually equivalent. If we choose an atlas of coordinate charts for which the metric is of Hessian form, we can induce a flat connection $D$ by differentiation with respect to the $u$-coordinates. The requirement that the transition maps be affine is exactly what is necessary for this connection to be globally defined.
Furthermore, we can induce the dual connection $D^*$ by differentiation with respect to the so-called dual coordinates $\theta$, which are defined as   
\begin{equation} \label{dual coordinates}
\theta_i := \frac{ \partial \Psi}{\partial {u^i}} .
\end{equation}
 In the dual coordinates, the metric is also of Hessian form, where the potential is the Legendre dual $\Psi^\ast$.
 For further details on this correspondence, we refer the reader to the book by Shima \cite{GHS}. 

There are topological and geometric obstructions for a given Riemannian manifold to admit a Hessian structure. In dimensions 4 and higher, there are local curvature obstructions as well (see \cite{AA}). As all of the manifolds of interest in this paper are open sets in $\mathbb{R}^n$ (which admit a global coordinate chart), we can construct Hessian metrics simply by choosing a convex potential.


\subsection{The Curvature of the Sasaki Metric}

Our primary interest is in the curvature of \textit{K\"ahler} Sasaki metrics. However, the most straightforward approach to calculating the curvature is to use the curvature formulas for a general Sasaki metric, and them simplify them using the dually flat structure.

For a Riemannian manifold $ (\M,g)$ with an affine connection $D$, it is of interest to understand the curvature of its associated Sasaki metric. Satoh calculated the full curvature tensor explicitly (Theorem 2.3 of \cite{HS}). For brevity, we will not present the full expression here.
However, in the case where $D$ is a flat connection, the formulas simplify considerably.

\begin{proposition}\label{Flat Sasaki metric}

Let $(\M, g ,D)$ be an affine manifold with flat connection $D$ and Levi-Civita connection $\nabla$. Let $\widetilde{R}$ be the Riemannian curvature tensor of the Sasaki metric $g ^D$ on $T \M$. For vector fields $\mathcal{X} , \mathcal{Y} , \mathcal{Z} , \mathcal{W} , \xi \in T_u \M$,
\begin{eqnarray*}
\widetilde {R}_{g^{D}}  \left( \mathcal{Z}_{\xi}^{H}, \mathcal{W}_{\xi}^{H}, X_{\xi}^{H}, \mathcal{Y}_{\xi}^{H} \right)=  {R}_{g^{D}}  \left( \mathcal{Z}_{\xi}^{V}, \mathcal{W}_{\xi}^{V}, \mathcal{X}_{\xi}^{V}, \mathcal{Y}_{\xi}^{V} \right) = &  R_{g}^{\nabla}( \mathcal{Z}, \mathcal{W},\mathcal{X}, \mathcal{Y}) 
\end{eqnarray*}
\begin{eqnarray*}
 \widetilde { R }_ {g ^ { D } } \left( \mathcal{Z} _ { \xi } ^ { H } , \mathcal{W} _ { \xi } ^ { V } , \mathcal{X} _ { \xi } ^ { V } , \mathcal{Y} _ { \xi } ^ { V } \right) = \widetilde { R } _ { g ^ { D } } \left( \mathcal{Z} _ { \xi } ^ { H } , \mathcal{W} _ { \xi } ^ { V } , \mathcal{X} _ { \xi } ^ { H } , \mathcal{Y} _ { \xi } ^ { H } \right) = 0
\end{eqnarray*}
\begin{eqnarray*}
\begin{aligned} \widetilde { R } _ { g } \left( \mathcal{Z} _ { \xi } ^ { H } , \mathcal{W} _ { \xi } ^ { V } , \mathcal{X} _ { \xi } ^ { H } , \mathcal{Y} _ { \xi } ^ { V } \right) = & - \frac { 1 } { 2 } \left( D _ { \mathcal{X} \mathcal{Z} } ^ { 2 } g \right) ( \mathcal{Y} , \mathcal{W} ) - \frac { 1 } { 2 } \left( D _ { \gamma ( \mathcal{X} , \mathcal{Z} ) } g \right) ( \mathcal{Y} , \mathcal{W} ) \\
 & + \frac { 1 } { 4 } \sum _ { i } \left( D _ { \mathcal{X} } g \right) \left( \mathcal{W} , e _ { i } \right) \cdot \left( D _ { \mathcal{Z} } g \right) \left( \mathcal{Y} , e _ { i } \right) \end{aligned}
\end{eqnarray*}

Here, $\{e_i \}$ is an orthonormal basis of $T_u \M$ and $\gamma ^ { D }$ is the difference between $D$ and the Levi-Civita connection on $\M$:
  $$\gamma ^ { D } ( \mathcal{X} , \mathcal{Y} ) = D _ { \mathcal{X} } \mathcal{Y} - \nabla _ { \mathcal{X} } \mathcal{Y}$$

\end{proposition}

When $(\M, g, D)$ is dually flat (i.e. Hessian), the situation simplifies further. To ease the computations, it is helpful to work in coordinates $\{ u^i \}_{i=1}^n$ where the Christoffel symbols of $D$ vanish. Doing so, we find the following identities.

\begin{enumerate}
\item The Riemannian metric $g$ is given by the Hessian of a convex potential $\Psi$.
\item In the induced coordinates $ \{ (u^i,v^i) \}_{i=1}^n$ on the tangent bundle $ (T \M,g^D,J^D)$, the complex structure can be written as
 $J^D \partial_{u^i} = \partial_{v^i}$ and $J^D \partial_{v^i} = - \partial_{u^i}$. As such, the coordinate chart $ (u^1,v^1, \ldots, u^n,v^n)$ is biholomorphic to an open set in $\mathbb{C}^n$ under the natural identification. 
\item There are simple expressions for the Riemannian curvature and the Christoffel symbols of the Levi-Civita connection.
\begin{enumerate}
\item The Riemannian curvature of the $(\M,g)$ is the following:
\[ R_g^\nabla (\partial_{u^i}, \partial_{u^j}, \partial_{u^k},\partial_{u^l}) = - \frac { 1 } { 4 }  \Psi ^ { p q } \left( \Psi _ { j l p } \Psi _ { i k q } - \Psi _ { i l p } \Psi _ { j k q } \right). \]
\item  The Christoffel symbols of the Levi-Civita connection satisfy the following identity:
\[ \Gamma_{ijk} = \frac{1}{2}\Psi_{ijk}  \hspace{1cm} \Gamma_{ji}^k = \frac{1}{2}\Psi_{ijm} \Psi^{km} . \]
\end{enumerate}
\item Using these formulas for the Christoffel symbols, we obtain a simple formula for $D_{\gamma^D( \mathcal{X},\mathcal{Z})}$ for two vectors $\mathcal{X}=\mathcal{X}^i \partial_{u^i}$ and $\mathcal{Z}=\mathcal{Z}^k \partial_{u^k}$.
\begin{eqnarray*}
D_{\gamma^D(X,Z)} 
			&=&  - \mathcal{X}^i \mathcal{Z}^k \Gamma_{ik}^r D_{\partial_{u^r}} = -\mathcal{X}^i \mathcal{Z}^k \Psi_{iks} \Psi^{sr} D_{\partial_{u^r}}
\end{eqnarray*}

\end{enumerate}

Combining these identities with the curvature formulas for the Sasaki metric, we find the following proposition.

\begin{proposition}[The curvature of a K\"ahler Sasaki metric] \label{Dually flat Sasaki metric}

Let $(\M,g,D)$ be a Hessian manifold. The Riemannian curvature of the Sasaki metric on $(T \M,g^D,J^D)$ is the following.
\begin{eqnarray*}
\widetilde {R}_{g^{D}}  \left( \partial_{u^i}\partial_{u^j},\partial_{u^k},\partial_{u^l} \right)=  \widetilde{R}_{g^{D}}   \left( \partial_{v^i},\partial_{v^j},\partial_{v^k},\partial_{v^l} \right) = &   - \dfrac { 1 } { 4 }  \Psi ^ { rs } \left( \Psi _ {jls } \Psi _ { iks } - \Psi _ { ilr } \Psi _ { jks } \right)
\end{eqnarray*}
\begin{eqnarray*}
 \widetilde { R } _ { g^D }  \left( \partial_{u^i},\partial_{v^j},\partial_{u^k},\partial_{v^l} \right) = & - \dfrac { 1 } { 2 } \Psi_{ijkl} + \dfrac { 1 } { 4 } \left( \Psi_{iks} \Psi^{sr} \Psi_{jlr}  \right) + \dfrac { 1 } { 4 }  \left( \Psi^{sr}  \Psi_{ jkr } \Psi_{ilr} \right) 
\end{eqnarray*}

Furthermore, when stated in terms of holomorphic vector fields, the curvature of $T \M$ satisfies the following identity:
\begin{eqnarray}\label{OHBC}
R_{i \bar j k \bar l}  & = &  \tilde R_{g^D}  \left( \partial_{u^i},\partial_{v^j},\partial_{u^k},\partial_{v^l} \right) - \tilde R_{ g^D} \left( \partial_{u^i},\partial_{u^j},\partial_{u^k},\partial_{u^l} \right) \\
  & =  & -\dfrac{1}{2} \Psi_{ijkl} + \dfrac{1}{2} \Psi^{rs}\Psi_{iks}\Psi_{jlr}. \nonumber
\end{eqnarray}

\end{proposition}

We remark that for a Hessian manifold, Shima defined the \text{Hessian curvature} to be the negative of formula \ref{OHBC} \cite{GHS}. We will not use this convention and instead work in terms of complex geometry.

\section{Optimal Transport and Complex Geometry}

With the background concluded, we can now state the central results of this paper, which relate the regularity theory of optimal transport to the complex geometry of the Sasaki metric.

\subsection{The MTW Tensor and the Curvature of $T \M$}
\hfill \\
When $c:X \times Y \to \mathbb{R}$ is a $\Psi$-cost, Ma, Trudinger and Wang observed that the MTW tensor takes the following form.
\begin{equation}\label{ConvexMTW}
 \mathfrak { S } _ { ( x , y ) } ( \xi , \eta )= (\Psi_{ijp}\Psi_{rsq}\Psi^{pq}-\Psi_{ijrs})\Psi^{rk} \Psi^{sl} \xi ^ { i } \xi ^ { j } \eta ^ { k } \eta ^ { l }. 
\end{equation}
In this formula, $k$ and $l$ are summed over, despite the double superscript resulting from the vector-covector ambiguity.

To make the connection between Equations \ref{OHBC} and \ref{ConvexMTW} precise, we induce $\M$ with the structure of a Hessian manifold. To do so, we use $\Psi$ as a potential for a Riemannian metric and let $D$ be the flat connection induced by differentiation with respect to the $u$-coordinates. This then induces $(T \M , g^D, J^D)$ with a K\"ahler metric, for which the following theorem is immediate.

\begin{theorem} \label{MTWisOAB}
Let $X$ and $Y$ be open sets in $\R^n$ and $c$ be a $\Psi$-cost. Then the MTW tensor satisfies the following identity:
\begin{equation} \label{MTWcurvature}
\mathfrak{S}(\xi, \eta) =  2 R_{g^D}(\xi, J \eta^\sharp , \xi, J \eta^\sharp) - 2 R_{g^D}(\eta^\sharp, \xi, \eta^\sharp, \xi). \end{equation}

Here, $R_{g^D}$ is the curvature  of the Sasaki metric on $(T \M, g^D, J^D)$ after sharping $\eta$ (recall that the $\eta$ in the MTW tensor is a covector with $\eta(\xi)=0$) and extending $\xi$ and $\eta$ to their real components in $T \M$. Furthermore, the cross curvature satisfies the same formula when we allow $\xi$ and $\eta$ to be arbitrary.
\end{theorem}

Note that it is important to be careful with the indices in the previous result. A previous version of this paper switched the roles of $j$ and $k$, which incorrectly showed a correspondence of the MTW tensor to the orthogonal bisectional curvature.

To simplify the discussion, we define the \textit{anti-bisectional curvature} by the following formula:
\begin{equation} \label{antibisectionalcurvature}
     \mathfrak{A} (\eta,\xi):= R_{g^D}(\xi, J \eta , \xi, J \eta) -  R_{g^D}(\eta, \xi, \eta, \xi)  \end{equation}
     
     In terms of holomorphic $(1,0)$-vectors $\xi$ and $\eta$, we have the following identity for the anti-bisectional curvature:
\[ \mathfrak{A}(\xi,\eta) = \frac{1}{2} \left( R_{g^D}( \xi, \overline \eta, \xi, \overline \eta) + R_{g^D} (\eta, \overline \xi, \eta, \overline \xi) \right). \]

Furthermore, we call the restriction of $\mathfrak{A}$ to orthogonal $\xi$ and $\eta$ the orthogonal anti-bisectional curvature, which we denote $\mathfrak{OA}$.
We say that the metric has (NOAB) if for all orthogonal real $\eta$ and $\xi$, $\mathfrak{OA}(\xi,\eta)$ is non-negative. We have chosen this terminology to reflect the fact that the bisectional curvature of two $J$-invariant planes satisfies the formula
\[ \mathfrak{H}(\eta,\xi) = R_{g^D}(\eta,  \xi, \eta, \xi) + R_{g^D}(\eta,  J \xi, \eta, J \xi),  \]
which differs from $\mathfrak{A}$ in that the sectional curvatures are added instead of subtracted. We conclude this discussion with the following remark, which we will use when we establish regularity for certain optimal transport problems.

\begin{remark}
The MTW tensor for a $\Psi$-cost is non-negative iff $T \M$ has (NOAB) on the set $T( X-Y) \subset T \M$. 
\end{remark}


\subsection{Relative $c$-Convexity of Sets and Dual Geodesic Convexity}
In order to establish regularity for optimal transport (as done in Theorem \ref{MTW}), not only is it necessary
to assume that the MTW tensor is non-negative, there are also assumptions about the
relative $c$-convexity of the supports of $\mu$ and $\nu$. For $\Psi$-costs, there is a natural geometric
interpretation for this notion, which we establish here.

Recall that for $x \in X$, a $c$-segment in $Y$ is the curve $c\textrm{-}\exp_x (\ell)$ for some line segment $\ell$ and a set $Y$ is c-convex relative to a set $X$ if, for all $x \in X$, $Y$ contains all $c$-segments between points in $Y$. For a $\Psi$-cost, relative $c$-convexity corresponds with geomdesic convexity with respect to the dual connection\footnote{Recall that the dual connection $D^*$ satisfies Equation  \ref{conjugateconnection}, where $D$ is the flat connection induced by differentiation with respect to the $u$ coordinates.} $D^*$.

To see this, we apply Equation \ref{dual coordinates}  to see 
\[ -c_{i,}= - \Psi_i (x-y) = - \theta_i(x-y), \]
where $\theta_i(x-y)$ is the point $x-y \in \M$ in terms of the dual coordinate $\theta_i$. By definition, $c$-segments correspond to straight lines in the $\theta$-coordinates, which are furthermore geodesics with respect to the dual connection $D^*$. This immediately implies the following.

\begin{proposition} \label{Relativecconvexity}
For a $\Psi$-cost, a set $Y$ is $c$-convex relative to $X$ if and only if, for all $x \in X$, the set $x-Y \subset \M$ is geodesically convex with respect to the dual connection $D^*.$ 
\end{proposition}

An analogous result holds for relative $c$-convexity of $X$ relative to $Y$. Combining the previous two results, we can restate Theorem \ref{MTW} in this new language.

\begin{theorem*}
 Suppose $X$ and $Y$ are smooth bounded domains in $\mathbb{R}^n$ and that $d \mu$ and $d \nu $ are smooth probability densities supported on $X$ and $Y$, respectively, bounded away from zero and infinity on their supports.
Consider a $\Psi$-cost for some convex function $\Psi: \M \to \mathbb{R}$ and suppose the following conditions hold.
\begin{enumerate}
\item $\Psi$ is $C^4$ and locally strongly convex (i.e. its Hessian is positive definite).
\item For all $x \in X$, $x-Y \subset \M$ is strictly geodesically convex with respect to the dual connection $D^*$. 
\item For all $y \in Y$, $X-y \subset \M$ is strictly geodesically convex with respect to the dual connection $D^*$.
\item The K\"ahler manifold $(T \M, g^D, J^D)$ has (NOAB) on the subset $T (X-Y)$.
\end{enumerate}

$\mathbb{T}_\textsc{u}$ be the $c$-optimal transport map carrying $\mu$ to $\nu$ as in Theorem 1.
Then $\textsc{u} \in C^\infty(\overline X)$ and $\mathbb{T}_\textsc{u}: \overline{X} \to \overline{Y}$ is a smooth diffeomorphism. 
\end{theorem*}

We should note that for many $\Psi$-costs of interest, $\Psi$ will not be uniformly strongly convex over its entire domain. This is no issue for the regularity theory, as we will restrict our attention to bounded sets $X$ and $Y$, so that $X-Y$ is precompact. As such, $\Psi$ will be strongly convex on $X-Y$.

\subsubsection{The Case of $\mathcal{D}_\Psi^{(\alpha)}$ Divergences}
\hspace{.01cm}

If we instead consider a divergence $\mathcal{D}_\Psi^{(\alpha)}:  \M \times \M  \to  \mathbb{R}$ (defined in \ref{Dalphadiv}), analogous results hold. In this setting, we construct a Hessian metric on $\M$ and consider $X$ and $Y$ as subsets of $\M$.  Theorem \ref{MTWisOAB} relates\footnote{There is a scaling factor of $\frac{1-\alpha^2}{2}$
between the curvature and the MTW tensor for a $\mathcal{D}_\Psi^{(\alpha)}$-divergence.} the MTW tensor of a  $\mathcal{D}_\Psi^{(\alpha)}$-divergence on $\M$ to the orthogonal anti-bisectional curvature of $T \M$ and Proposition \ref{Relativecconvexity}. From an information geometric point of view, this is a more natural construction as it eliminates the need to consider the difference set $X-Y$. Furthermore, when the potential is a log-partition function for an exponential family, its domain is guaranteed to be convex so the divergence is well-defined.

\section{Applications} 

As the results in the previous section give a new interpretation for prior work, it is natural to ask for new results that can be found using this approach. In this section, we give several such applications.  We will not provide the derivations of the identities in this section, as they are very involved. In order to compute the associated curvature tensors, we have written a Mathematica notebook, which is available online \cite{MNB}.

\subsection{A Complete, Complex Surface with (NOAB)}

One question of considerable interest in complex geometry is to understand complete K\"ahler metrics with various non-negativity properties. Most famously, the Frankel conjecture states that if a compact K\"ahler manifold has positive holomorphic bisectional curvature, then it is biholomorphic to the complex projective space $\mathbb{CP}^n$. This conjecture was independently proven by Mori \cite{Mori} in 1979 and Siu-Yau \cite{SY} 1980.

 For compact K\"ahler manifolds, it is possible obtain this result under weaker curvature assumptions. For instance, all compact K\"ahler manifolds with positive orthogonal bisectional curvature are biholomorphic to $\mathbb{CP}^n$ (see \cite{XXC} \cite{FLW} \cite{GZ}). Furthermore, all compact irreducible K\"ahler manifold with non-negative isotropic curvature are either Hermitian symmetric or else biholomorphic to $\mathbb{CP}^n$ \cite{HSesh}. For complex surfaces, non-negative orthogonal bisectional curvature is equivalent to non-negative isotropic curvature\footnote{In higher dimensions, non-negative isotropic curvature is a stronger assumption than non-negative orthogonal bisectional curvature.} \cite{LN} so the previous result gives a classification of such surfaces. For non-compact manifolds, it is natural to ask whether similar results hold. The most famous conjecture in this direction is Yau's uniformization conjecture, which states that any complete irreducible non-compact K\"ahler metric of non-negative bisectional curvature is biholomorphic to $\mathbb{C}^n$ \cite{STY}. Although the full conjecture is still open, Liu proved it under certain volume growth assumptions \cite{GLY}. 
 
Since the anti-bisectional curvature appears similar to the bisectional curvature, we can also ask how much control non-negative anti-bisectional curvature or (NOAB) provides over the geometry of a K\"ahler manifold. Using the polarization formula \cite{Hawley}, it can be shown that any metric of constant holomorphic sectional curvature has vanishing orthogonal anti-bisectional curvature, so any Hermitian symmetric space satisfies (NOAB). One can then ask whether there are other examples.

The following example gives a very interesting metric which satisfies (NOAB) and completeness but is neither Hermitian symmetric nor biholomorphic to $\mathbb{C}^n$.

\begin{example}[A complete surface with (NOAB)] \label{NOBexample}
Consider the negative half-plane $\M =\mathbb{H} := \{ (u_1, u_2) ~|~ u_2 < 0 \} $. Induce it with a Hessian metric with the potential function $\Psi: \mathbb{H} \to \mathbb{R}$ given by
\[  \Psi(u)= - \frac{u_1^2}{4u_2}-\frac{1}{2} \log(-2 u_2). \] 
\end{example}

 For a vector $\xi =\partial_{u_1}+ a \partial_{u_2}$ and a covector $\eta=a du_1-du_2$, the associated orthogonal anti-bisectional curvature\footnote{By computing anti-bisectional curvature solely on real vectors and covectors, we are slightly abusing notation. To formalize this, extend $\xi$ and $\eta$ to their real counterparts on $T \M$.} on $T \mathbb{H}$ is given by
\[ \mathfrak{OA}(\eta^\sharp,\xi) = \frac{6 a^2( -a u_1^2 +u_2)^2}{u_2^2}. \]
As such, the metric has (NOAB).
This metric is of independent interest, and for a more complete discussion, we refer the reader to the work of Molitor \cite{Molitor}. We will note a few of its curvature properties in passing.
For a vector $\xi=\partial_{u_1}+ a \partial_{u_2}$ and a covector $\eta= du_1+ a du_2$, the holomorphic sectional curvature is given by
\[ \mathfrak{H}(\eta,\xi) =2-4a^2+a^2 \left( -8 + \frac{6u_1^2}{u_2^2} \right) - 12 \frac{a u_1}{u_2}. \]
As such, the holomorphic sectional curvature does not have a definite sign. It can similarly be shown that the orthogonal bisectional curvature also does not have a definite sign. However, the metric does have constant negative scalar curvature. This manifold is complete and Stein (it is biholomorphic to an open set in $\mathbb{C}^2$). However, it has the standard complex structure on a half-space in $\mathbb{R}^4$, so is \textit{not} biholomorphic to $\mathbb{C}^2$.

\subsubsection{The Fisher metric of the normal distribution}

Although this example has interesting theoretical properties, it may appear to be a somewhat ad hoc construction without context. In fact, it is a natural example from information geometry. If we consider $u_1$ and $u_2$ as the natural parameters of the normal distribution (i.e. $u_1= \dfrac{\mu}{\sigma^2}$  and  $u_2= \dfrac{-1}{2\sigma^2}$), then the Riemannian metric $ g_{ij} = \Psi_{ij} $ 
 is the Fisher metric on the statistical manifold of univariate normal distributions.
 As a Riemannian manifold, $(\mathbb{H}, g)$ is a complete hyperbolic surface (which motivated our choice of notation). Note, however, that the $(u_1,u_2)$ coordinates do \textit{not} induce the standard half-plane model of hyperbolic space.\footnote{In $(\mu,\sigma)$ coordinates, the Fisher metric is $ds^2= \dfrac{1}{\sigma^2} \left( d \mu^2 + 2 d \sigma^2 \right)$, which is much closer to the standard model. }

\subsubsection{A closely related example}

Using the normal statistical manifold, it is possible to construct another K\"ahler metric which satisfies (MTW). This space is actually Hermitian symmetric and was first constructed by Shima (Example 6.7 of \cite{GHS}).

Consider the domain 
\[ \widetilde{ \mathcal{M}} := \{ (\theta_1,\theta_2) ~|~ \theta_2-\theta_1^2>0 \} \]
and induce it with a Hessian metric with potential $\Psi^*(\theta) =-\frac{1}{2} -\log(\theta_2-\theta_1^2)$.


This potential arises from the parameterization of the univariate normal distribution in terms of its dual parameters $\theta_1 = \mu$ and $\theta_2 = \mu^2 + \sigma^2$ and the potential $\Psi^*(\theta)$ is the Legendre dual of the above potential in Example \ref{NOBexample}. Computing the anti-bisectional curvature for a vector $\xi$ and covector $\eta$, we find that it satisfies
\[ \mathfrak{A}(\xi,\eta^\sharp) = - \eta(\xi)^2. \] 

As such, the orthogonal anti-bisectional curvature vanishes and the holomorphic sectional curvature is a negative constant. From this, we can see that the geometry of $T \widetilde{ \mathcal{M}}$ is of independent interest, as it is a complete Hermitian symmetric space with constant negative holomorphic sectional curvature. This is an example of a Siegel upper half-space, but we will postpone a more complete discussion to future work.


We note that it is possible to construct other K\"ahler metrics with (NOAB) that are very similar to $T \widetilde{ \mathcal{M}}$. Using a similar construction for round multivariate Gaussian distributions, it is possible to construct such a Hermitian symmetric space in arbitrary dimensions. For another example, we can consider the potential $\Psi(\theta_1,\theta_2) =-\frac{1}{2} -\log(\theta_2-\theta_1^4)$, which also has (NOAB).

\subsubsection{Regularity for an associated cost function}

We can also use this potential to construct a cost function with a natural regularity theory. Instead of using a $\Psi$-cost, it is more natural to consider the $\mathcal{D}^{(\alpha)}_\Psi$-divergence. For this cost function, we can apply our previous calculations to obtain the following result.

\begin{cor*}
Suppose $\mu$ and $\nu$ are probability measures supported on bounded subsets $X$ and $Y$ of the normal statistical manifold $\M$. Suppose further that the following regularity assumptions hold.
\begin{enumerate}
\item $\mu$ and $\nu$ are absolutely continuous with respect to the Lebesgue measure. Furthermore, $d \mu$ and $d \nu$ are smooth and bounded away from zero and infinity on their respective supports.
\item For all $x \in X$, $\frac{x+ Y}{2}$ is strictly convex with respect to the coordinates $\theta_1 = \mu$ and $\theta_2 = \mu^2+\sigma^2$. Furthermore, the same property holds for $\frac{X+ y}{2}$ for all $y \in Y$.
\end{enumerate}
 Let $ c (x,y)$ be the cost function given by
\[ c(x,y)=  \frac{1}{2} \Psi(x) +  \frac{1}{2} \Psi(y) - \Psi \left( \frac{x+y}{2} \right), \]
where $\Psi$ is the convex function given in Example \ref{NOBexample}.
Then the $c$-optimal map $\mathbb{T}_\textsc{u}$ taking $\mu$ to $\nu$ is smooth.
\end{cor*}

\subsection{The Regularity of Pseudo-Arbitrages}

Recently, a series of papers by Pal and Wong (\cite{ECF}-\cite{GRA}, \cite{POT}-\cite{LD}) have studied the problem of finding \textit{pseudo-arbitrages}, which are investment strategies which outperform the market portfolio under ``mild and realistic assumptions." Their work combines information geometry with optimal transport and mathematical finance to reduce the problem to solving optimal transport problems where the cost function is given by a so-called log-divergences.

A central result in \cite{ECF} shows that a portfolio map $\mathbf{ \pi}$ outperforms the market portfolio almost surely in the long run iff it is a solution to the Monge problem for the cost function $ c: \R^{n-1} \times \R^{n-1} \to \R$ given by
\begin{equation}  \label{ECF cost}
   c (x, y) := \log \left( 1 + \sum_{i=1}^{n-1} e^{x^i -y^i} \right) - \log(n) - \frac{1}{n} \sum_{i=1}^{n-1} x^i - y^i. 
   \end{equation}
   
To give some context for this cost function, it is instructive to consider $x$ and $y$ as the natural parameters of the multinomial distribution. For natural parameters $\{ x^i \}_{i=1}^{n-1}$, we can compute the probability $p_i$ of the $i$-th event (in this context, the $i$-th market weight) using the following formulas:
\begin{eqnarray*}
 p_i &= & \frac{e^{x^i}}{1+ \sum_{j=1}^{n-1} e^{x^j}}\textrm{ for } 1 \leq  i<n, \\
 p_n & = & \frac{1}{1+ \sum_{j=1}^{n-1} e^{x^j}}. 
\end{eqnarray*}

To write this cost function in a more familiar form, we similarly find probabilities $q_i$ associated to the $y$-parameters and fix $\boldsymbol \pi = \left(\frac{1}{n}, \ldots \frac{1}{n} \right) \in \triangle^n$. Rewriting our cost in these terms, we have the following:
 \[ T(\,p\,|\,q\,) = \log \left( \sum_{i=1}^n \boldsymbol \pi_i \frac{p_i}{q_i} \right) - \sum_{i=1}^n  \boldsymbol \pi_i \log \left( \frac{p_i}{q_i} \right) \]

This quantity is known as the free energy in statistical physics \cite{ECF} and by various different names in finance (such as the ``diversification return", the ``excess growth rate," the ``rebalancing premium" and the ``volatility return").
 Since Pal and Wong refer to this as a \textit{logarithmic divergence}, we refer to this cost as the logarithmic cost. 
This cost function is not symmetric, so is not induced by any distance function. However, Jensen's inequality shows that it is a divergence.

The main focus of Pal and Wong's work is to study the information geometric properties of divergences functions induced by log-convex functions, of which $T(\, \cdot \, | \, \cdot \,)$ is only a single example. For any log-convex function, one can define a corresponding divergence $D[ \, \cdot \, | \, \cdot \, ]$, which has a self-dual representation in terms of the logarithmic cost (see Proposition 3.7 of \cite{ECF}). 
 In order to study optimal transport, we do not specify the log-convex function a priori. In fact, such a function induces the {\em solution} to an optimal transport problem.

For the logarithmic cost, only the first term affects optimal transport. As such, we instead consider the cost function 
\[  \tilde c (x, y) := \log \left( 1 + \sum_{i=1}^{n-1} e^{x^i -y^i} \right).  \]

This is now a $\Psi$-cost for the convex function 
\[ \Psi (u) =  \log \left( 1 + \sum_{i=1}^{n-1} e^{u^i } \right). \]

As such, we can apply Theorem \ref{MTWisOAB} to compute the MTW tensor for the cost $\tilde c$. For a vector $\xi$ and a covector $\eta$, the anti-bisectional curvature of $T \R^{n-1}$ (denoted $ \mathfrak{A}$) with Hessian metric induced by $\Psi$ is
 \[ \mathfrak{A}(\xi,\eta^\sharp) = 2 g(\eta^\sharp,\xi)^2. \]

As such, the MTW tensor identically vanishes and the cost has non-negative cross-curvature. A proof for this identity can be found in Proposition 3.9 of Shima's book \cite{GHS}. From the curvature formulas, we see that this potential induces a K\"ahler metric on $\mathbb{C}^n$ with constant positive holomorphic sectional curvature. As such, it is an Hermitian symmetric space, but is not complete.

In order to apply the Trudinger-Wang result, we must also determine what relative $c$-convexity means in this context. To do so, we solve for the dual coordinates to the natural parameters $u^i$ by calculating $\partial_{u^i} \Psi$ for $i=1, \ldots, n-1$. Doing so, we find that the dual coordinates are
\begin{eqnarray*}
	 \frac{e^{x^i}}{1+ \sum_{j=1}^{n-1} e^{x^j}}\textrm{ for } 1 \leq  i<n, 
\end{eqnarray*}
which are exactly the formulas for the market weights $p_i$. This is initially surprising, but has a natural interpretation in terms of information geometry

\subsubsection{The information geometry of the multinomial distribution}
It is worth discussing the geometry of this example in more detail. It turns out that if we consider the $\{x^i\}$ as natural parameters, the potential $\Psi$ induces the Fisher metric of the multinomial distribution, which is a standard example in statistics. Geometrically, this is the round metric on the positive orthant of a sphere, which immediately shows that neither the underlying Hessian metric nor the Sasaki metric is complete.

It is worth mentioning that this metric cannot be extended to a K\"ahler Sasaki metric on the tangent bundle of the entire sphere, due to the fact that the sphere is \textit{not} an affine manifold. 

 For an exponential family of probability distributions, the dual coordinates are the expected values of the natural sufficient statistics. More specifically, for the multinomial distribution the dual coordinates are precisely the original market weights, which explains the relationship between the market weights and the partial derivatives of the potential function. As such, if we let $\mathcal{P}$ be the coordinate transformation from the natural parameters $x$ to the market weights $p$, a subset $X \subset \R^{n-1}$ is relatively $c$-convex iff the set $ \mathcal{ P} (X) $ is convex as a subset of the probability simplex in the usual sense. Using this transformation, we say that a subset $ \mathcal{P} (X) $ of the probability simplex has \textit{uniform probability} if $X$ is a precompact set. More concretely, a subset  $\mathcal{P} (X) $ has uniform probability if and only if there exists $\delta>0$ so that for all $p \in \mathcal{P} (X) $ and $1 \leq i \leq n$, $p_i > \delta.$

 \subsubsection{Regularity of optimal transport}

 From these observations and the previous identity for the MTW tensor, we can derive the following regularity result.

\begin{cor} \label{Simplexregularity}
Suppose $\mu$ and $\nu$ are smooth probability measures supported respectively on subsets $X$ and $Y$ of the probability simplex $\triangle$. Suppose further that the following regularity assumptions hold:
\begin{enumerate}
\item $X$ and $Y$ are smooth and strictly convex. Furthermore, both have uniform probability (as defined above).
\item $\mu$ and $\nu$ are absolutely continuous with respect to the Lebesgue measure and $d \mu$ and $d \nu$ are bounded away from zero and infinity on their supports.
\end{enumerate}
 Let $\hat c (p,q)$ be the cost function given by
\[ \hat c(p,q) = \log \left( \frac{1}{n} \sum_{i=1}^n \frac{q_i}{p_i} \right) - \frac{1}{n} \sum_{i=1}^n \log \frac{q_i}{p_i}. \]
Then the $\hat c$-optimal map $\mathbb{T}_\textsc{u}$ taking $\mu$ to $\nu$ is smooth.
\end{cor}

 In \cite{MSP}, Pal and Wong study the cost function $\hat c$ and use it to define a displacement interpolation between two probability measures. In their paper, they inquire about the regularity problem for this interpolation. We can now answer this question using the previous result.

\begin{cor}
Suppose $\mu$ and $\nu$ are smooth probability measures satisfying the assumptions of Corollary \ref{Simplexregularity} and that $\mathbb{T}_\textsc{u}$ is the $\hat c$-optimal map transporting $\mu$ to $\nu$. Suppose further that $\mathbb{T}(t) \mu$ is the displacement interpolation from $\mu$ to $\nu$ defined by $\mathbb{T}(t) = t \cdot \textrm{Id}+ (1-t) \mathbb{T}_\textsc{u}$. Then $\mathbb{T}(t)$ is smooth, both as a map on the probability simplex for fixed $t$ and also in terms of the $t$ parameter.
\end{cor}

 For $t=1$, the solution to the interpolation problem is simply $\mathbb{T}_\textsc{u}$ and so Corollary \ref{Simplexregularity} shows that the potential \textsc{u} is smooth. Since the displacement interpolation defined in \cite{MSP} linearly interpolates the potential functions, the associated displacement interpolation is also smooth for $0 \leq t \leq 1$. 

In closing, we note that the cost function considered here is very similar, but not identical, to the radial antennae cost, which was studied by Wang \cite{RAP}. It is of interest to determine whether there is some deeper connection between these two costs which explains their apparent similarity.


\subsection{Other Examples in Complex Geometry and Optimal Transport}

While writing this paper, we were able to find several more examples of Hessian manifolds whose tangent bundles have non-negative bisectional curvature or (NOAB).

 Relatively few examples of positively curved metrics are known (for some examples, see \cite{WZ}), so this method may be helpful for finding new ones. One limitation of this approach is that many of the manifolds are not complete as metric spaces. It would be of interest to determine which convex functions induce complete K\"ahler metrics with non-negative or positive orthogonal bisectional curvature, and we plan to study this problem in future work.

Each of these examples further induces a cost with non-negative MTW tensor. Furthermore, since many of these examples are obtained from statistical manifolds, it may be possible to use them to induce meaningful statistical divergences.

\begin{enumerate}

\item $ \Psi(u) = -\log \left( 1-\sum_{i=1}^n e^{u_i} \right)$. This potential induces a Sasaki metric with constant negative holomorphic sectional curvature.
\[ \mathfrak{A}(\xi,\eta^\sharp) = - \eta(\xi)^2. \]

From an information geometric point of view, this is the Fisher metric of the negative multinomial distribution. As a Hessian manifold, $\M$ is a non-compact metric of constant negative holomorphic sectional curvature, but is not complete.

\item $\Psi(u)=\left( e^{u_1}+e^{u_2} \right)^p$  for $0<p<1$. For a vector $\xi = \xi_1 \partial_{u_1}+ \xi_2 \partial_{u_2}$ and covector $\eta= \eta_1 d{u_1}+\eta_2 d{u_2}$, the associated orthogonal anti-bisectional curvature of the Sasaki metric is given by
\[ \mathfrak{OA}(\xi,\eta^\sharp)= \frac{2(1/p-1)(a-1)^2 (e^{u_1}+ae^{u_2})^2}{(e^{u_1}+e^{u_2})^{2+p}}. \]

As Hessian manifolds, these metrics are non-compact, but are not complete.

\item $ \Psi(u)= \log( \cosh(u_1)+\cosh(u_2))$. This potential induces a Sasaki metric whose bisectional curvature is non-negative.
For a vector $\xi = \xi_1 \partial_{u_1}+ \xi_2 \partial_{u_2}$ and covector $\eta= \eta_1 d{u_1}+\eta_2 d{u_2}$, the associated bisectional curvature is given by
\[ \mathfrak{H}(\xi,\eta^\sharp)= |\xi |^2| \eta|^2+4 \xi_1 \xi_2  \eta_1 \eta_2 \]
where $|\xi|^2 = \xi_1^2+ \xi_2^2$ and $|\eta|^2 = \eta_1^2+ \eta_2^2$.

Furthermore, the anti-bisectional curvature also satisfies the same formula.
\[ \mathfrak{A}(\xi,\eta^\sharp)= |\xi |^2| \eta|^2+4 \xi_1 \xi_2  \eta_1 \eta_2. \]

As such, this metric has (NOAB) and non-negative bisectional curvature. As a Hessian manifold, this metric is bounded, and so is not complete. Note that the curvature of this metric is in fact parallel with respect to $D$, which makes it an interesting example. We will explore this metric further in future work.

\end{enumerate}

\section{Open questions}

\subsection{The Complex Geometry of Optimal Maps}

It is of interest to determine whether the regularity theory for optimal transport with $\Psi$-costs can be established in terms of complex Monge-Ampere equations. To do so, we would first need to find a version of Theorem 1 for complex Monge-Ampere equations, which we leave for future work.
For a more complete overview on complex Monge-Ampere equations, we refer the reader to the paper by Phong, Song and Sturm \cite{PSS}.

One striking feature of the pseudo-Riemannian theory of optimal transport is the natural geometric interpretation for optimal maps. More precisely, if one deforms the pseudo-metric by a particular conformal factor (which is determined by the respective densities), the optimal map is induced by a maximal codimension-one surface with respect to the conformal pseudo-metric \cite{KMW}. At present, we do not have a corresponding interpretation for this in the complex setting, and we leave this question for future work.

\subsection{Implications for Optimal Transport}

The primary focus of this work is to use optimal transport theory to study complex geometry and information geometry. However, it remains an open question what can be proven about optimal transport using this approach. For instance, it is hard to find cost functions which satisfy $MTW(0)$. Similarly, relatively few K\"ahler metrics of positive curvature are known and there are various stability and gap theorems about them (see, e.g., \cite{GLY} and \cite{NN}). In optimal transport, we expect that similar results can be proven using complex geometry. We plan to explore this topic in future work.

\subsection{A Potential Non-K\"ahler Generalization}

While $\Psi$-costs yield many interesting examples, there are many relevant cost functions which are not of this form. As such, one natural generalization of the construction considered here is to instead consider a Lie group $\mathcal{G}$ and cost functions of the form $\Psi(x\cdot y^{-1})$ for $x,y \in \mathcal{G}$. Our work thus far can be interpreted as doing this calculation in the special case where $\mathcal{G}$ is Abelian. For non-Abelian groups, we hope it is possible to recover the MTW tensor as a curvature tensor in almost complex geometry. In this case, there would be correction terms due to the non-Abelian nature of the group and the almost complex structure on the tangent bundle $T \mathcal{G}$ would fail to be integrable, so the associated theory would be non-K\"ahler.

The motivation for this approach comes from the $2$-Wasserstein cost. For this cost, many of the known examples satisfying $MTW(\kappa)$ are compact Lie groups. For instance, it is known that $SO(3) \cong \mathbb{RP}^3$ satisfies $MTW(\kappa)$. We hope that this approach can be used to show the MTW property for other compact Lie groups.

\subsection{The Complex Geometry of the Anti-Bisectional Curvature}

The orthogonal anti-bisectional curvature is a very subtle invariant, and its geometry is quite mysterious. It does not determine the full curvature tensor of the manifold, since all spaces with constant holomorphic sectional curvature have vanishing orthogonal anti-bisectional curvature. However, it is non-negative for some important metrics which don't have any other obvious non-negativity properties. In future work, we hope to understand this curvature more fully and to understand what sort of control it exerts over the geometry of a complex manifold.

\end{document}